\documentclass[10pt]{article}
\usepackage{amssymb}
\usepackage{latexsym}
\usepackage{amsthm}
\usepackage{amscd}
\usepackage{amsmath}
\usepackage{mathrsfs}
\usepackage[all]{xy}
\usepackage{graphicx}
\RequirePackage{amsthm}

\newtheorem{prop}{\bfseries Proposition}
\newtheorem{theo}[prop]{\bfseries Theorem}
\newtheorem{lemm}[prop]{\bfseries Lemma}
\newtheorem{coro}[prop]{\bfseries Corollary}
\theoremstyle{definition}

\newtheorem{rema}[prop]{Remark}
\newtheorem{exem}[prop]{Example}
\numberwithin{prop}{section}

\def\ZZ{{\bf Z}}

\def\Ca{{{\mathscr C}}}
\def\A{{{\bf A}^1_k}}
\def\P{{{\bf P}^1_k}}
\def\OO{{\mathscr O_X}}
\def\oX{{\mathscr L}}
\def\M{{\mathscr B}}
\def\E{{\mathscr F}}
\def\G{{\mathscr G}}

\def\MM{{\mathscr M}}

\def\W{{\mathscr W}}
\def\Z{{\mathscr Z}}
\def\H{{\mathscr H}}
\def\E{{\mathscr E}}
\def\hyp{{\rm hyp}}

\DeclareMathOperator{\coker}{coker}
\DeclareMathOperator{\Div}{div}
\DeclareMathOperator{\res}{Res}
\DeclareMathOperator{\proj}{Proj}

\DeclareMathOperator{\DDiv}{Div}
\DeclareMathOperator{\id}{id}

\DeclareMathOperator{\im}{Im}

\begin{document}
\title{On the $1$-pointed curves arising as \'etale covers of the affine line in positive characteristic}
\author{Leonardo Zapponi\\ EPFL, Lausanne}
\maketitle

\section*{Introduction}

 Let $k$ be an algebraically closed field of positive characteristic. The goal of this paper is to characterize the proper smooth curves $X/k$ of positive genus $g$ equipped with a $k$-rational point $P$ such that $X\setminus P$ can be realized as an \'etale cover of the affine line. It turns out that this question is closely related to the action of the Cartier operator on the set of differential forms on $X$. A first analysis shows that such a cover exists if and only if there exists an exact regular differential form $\omega$ on $X$ having a unique zero at $P$ (of order $2g-2$). The essential  inconvenient of this (elementary) approach is that it doesn't give any control of the degree of the cover. A finer investigation, essentially based on the duality between the Cartier operator and the Frobenius map, allows us to overcome this problem. The final characterization is particularly simple and adapted for explicit calculations.

The paper is organized as follows: in the first section we fix the notation and recall some basic facts concerning invertible sheaves on a curve. In the second section we start the investigation by some elementary methods, only using the Riemann-Roch theorem. In the third section we work in a more sheaf-theoretic setting and obtain some duality resutls which are used in section 4 when studying the question of the minimal degree of the cover. In the last section we give many explicit examples and complements which illustrate the general results obtained.

\section{Some notation} 

Let $k$ be an algebraically closed field of positive characteristic $p$. Throughout this paper,  we fix an irreducible, proper and smooth curve $X/k$ of positive genus $g$. We denote by $\OO$ its structure sheaf and by $K$ its field of rational functions. For simplicity, given a sheaf of abelian groups $\mathscr F$ on $X$, we write $H^i(\mathscr F)$ instead of $H^i(X,\mathscr F)$. For any divisor $D$ on $X$, let $\oX(D)$ be the invertible sheaf whose sections over an open subset $U\subset X$ are the rational functions $x\in K$ such that $v_P(x)\geq-v_P(D)$ for any $P\in U$. As usual, let $L(D)$ be the set of its global sections and denote by $l(D)$ its dimension, as a $k$-vector space. Let $\Omega^1_X$ the sheaf of regular K\"ahler differentials on $X$. According to the notation of \cite{Serre_1}, let $\Omega(D)$ be the set of global sections of the sheaf $\Omega^1_X(D)=\Omega^1_X\otimes\oX(-D)$ and denote by $i(D)$ its dimension so that the Riemann-Roch formula can be written as $l(D)-i(D)=\deg(D)-g+1$.

Consider the absolute Frobenius morphism $F:X\to X$. It induces a homomorphism of $\OO$-modules $\OO\to F_*\OO$ and thus a semi-linear map $H^1(\OO)\to H^1(\OO)$ (since the Frobenius morphism is affine, we have a natural isomorphism $H^i(F_*\mathscr F)\cong H^i(\mathscr F)$ for any coherent sheaf $\mathscr F$ on $X$). Recall that the {\it Cartier operator} is a homomorphism of $\OO$-modules $\Ca:F_*\Omega^1_X\to\Omega^1_X$ inducing a semi-linear map $\Ca:\Omega^1_K\to\Omega^1_K$ (where we have set $\Omega^1_K=\Omega^1_X\otimes  K$) satisfying the two following properties (cf.~\cite{Cartier}, \cite{Colliot} or \cite{Yui}): given a meromorphic section $\omega\in\Omega^1_K$, we have
\begin{itemize}
\item $\Ca(\omega)=0$ if and only if $\omega$ is exact, i.e. $\omega=dx$ with $x\in K$.
\item $\Ca(\omega)=\omega$ if and only if $\omega$ is logarithmic, i.e. $\omega=dx/x$ with $x\in K$.
\end{itemize}
As usual, the dual of a $k$-vector space $V$ is denoted by $V^\vee$. If we identify $\Omega(0)$ with $H^1(\OO)^\vee$ by Serre duality then the Cartier operator is just the adjoint of the Frobenius map (cf.~\cite{Serre_2}). More precisely, there is a natural perfect pairing
$$\aligned H^1(\OO)\times\Omega(0)\to k\\
(r,\omega)\mapsto\langle r,\omega\rangle
\endaligned$$
which satisfies the relation
$$\langle r^p,\omega\rangle=\langle r,\Ca(\omega)\rangle^p$$

\section{Etale covers of the affine line. Elementary approach}

We start with some general considerations concerning the exact differential forms on $X$. All the results of this section follow from the Riemann-Roch theorem. No more sophisticated  techniques are needed. First of all, for any meromorphic differential form $\omega\in\Omega^1_K$, denote by ${\rm Supp}(\omega)\subset X $ the support of ${\rm div}(\omega)$, i.e. ${\rm Supp}(\omega)={\rm Supp}(\omega)_0\cup{\rm Supp}(\omega)_\infty$, where ${\rm Supp}(\omega)_0$ (resp. ${\rm Supp}(\omega)_\infty$) is the set of zeroes (resp. of poles) of $\omega$.

\begin{lemm}\label{lemm1} Let $\omega$ be an exact meromorphic differential form on $X$. For any point $P\in X $ there exists a rational function $x\in K$ which is regular outside ${\rm Supp}(\omega)_\infty\cup P $ and satisfies $dx=\omega$. We can moreover require the condition
$$v_P(x)\geq\min\{(1-2g)p,v_P(\omega)+1\}$$ 
\end{lemm}

\begin{proof} Let $y\in K$ be a rational function such that $dy=\omega$. Consider an element $Q\in X $ not belonging to ${\rm Supp}(\omega)_\infty\cup P $ and suppose that $y$ has a pole at $Q$ of order $n>0$. Then, $n$ is a multiple of $p$ (the contrary would imply that $dy$ has a pole at $Q$ of order $n+1$). Setting $n=mp$, the Riemann-Roch theorem asserts that there exists a rational function $t\in L((2g-1)P+mQ)$ having a pole at $Q$ of exact order $m$. Indeed, the divisor $D=(2g-1)P+mQ$ has degree greater than or equal to $2g$ and thus $i(D)=i(D-Q)=0$. This gives $l(D)-l(D-Q)=1$, so that there exists a rational function $t$ satisfying the above properties. In particular, we can find an element $c\in k^*$ such that the rational function $y_1=y-ct^p$ has a pole at $Q$ of order strictly less than $n$. Remark that $dy_1=dy=\omega$. By iteration of the above procedure, we obtain a rational function $y'$ which is regular at $Q$ and satisfies $dy'=\omega$. Proceeding in this way for all the remaining poles of $y'$ not belonging to ${\rm Supp}(\omega)_\infty\cup P $, we finally get a rational funtion $x$ which is regular outside ${\rm Supp}(\omega)_\infty\cup P $. We clearly have $v_P(x)\leq v_P(\omega)+1$. Suppose finally that $v_P(x)<\min\{(1-2g)p,v_P(\omega)+1\}$. Since $dx=\omega$, we obtain $v_P(x)=np$ with $n\leq-2g$. In particular, there exists a rational function $t$ having a unique pole at $P$ of order $n$ such that $v_P(x-t^p)>v_P(x)$. If we iterate this construction we obtain a rational function satisfying the conditions of the lemma.
\end{proof}

\begin{prop}\label{prop1}For any nonempty finite subset $S\subset X $, the following conditions are equivalent:

\begin{enumerate}
\item There exists an \'etale cover $X\setminus S\rightarrow\A$.
\item There exists an exact meromorphic differential form $\omega$ whose support is contained in $S$.
\end{enumerate}
\end{prop}

\begin{proof}
Let's start by supposing that there exists an \'etale cover  $X\setminus S\rightarrow\A$. It can be uniquely extended to a ramified cover $\beta:X\to\P$. Set $\P=\A\cup\infty$. Consider a non-constant function $x\in L(\infty)$ (the Riemann-Roch theorem gives $l(\infty)=2$). The differential form $dx$ has a unique pole at $\infty$ and is regular and nowhere vanishing on $\A$. In particular, its pull-back $\omega=\beta^*dx\in\Omega^1_K$ is regular and nowhere vanishing on $X\setminus S$, so that ${\rm Supp}(\omega)\subset S$. Moreover, we clearly have $\Ca(\omega)=\Ca(\beta^*dx)=\Ca(d\beta^*x)=0$ and thus $\omega$ is exact. Conversely, suppose that the second condition is fullfilled. Then, fixing an element $P\in S$, Lemma~\ref{lemm1} implies that there exists a rational function $y$ which is regular outside ${\rm Supp}(\omega)_\infty\cup P \subset S$ and satisfies $dy=\omega$. Let $t\in K$ be a rational function having $S$ as set of poles (it exists by the Riemann-Roch theorem and we can moreover suppose that the inequality $v_Q(t)=-1$ holds for any $Q\in S\setminus P $ and that $v_P(t)\geq 1-2g$) and set $x=y+ct^p$, with $c\in k^*$. By construction, there exists $c$ such that $S$ is the set of poles of $x$, so that we obtain a finite cover $X\setminus S\rightarrow\A$, which is \'etale, since the zeores of $dx=\omega$ are contained in $S$.
\end{proof}

\begin{coro}\label{cor1} There exists a finite set $S\subset X $ such that the curve $X\setminus S$ can be realized as an \'etale cover of the affine line.
\end{coro}

\begin{proof} Just take $x\in K\setminus K^p$, $\omega=dx$,
$S=\text{Supp}(\omega)$ and apply Proposition~\ref{prop1}.
\end{proof}

We now investigate the case where $S$ is reduced to a single element, which is the main interset of the paper. From now on, we restrict to this situation. 

\begin{theo}\label{th1} For any element $P\in X $, the following conditions are equivalent:
\begin{enumerate}
\item There exists an \'etale cover $X\setminus P \rightarrow\A$.
\item There exists an exact differential form which is regular and nowhere vanishing outside $P$.
\item The $k$-vector space $\Omega((2g-2)P)$ is non-trivial and $\Ca(\Omega((2g-2)P))=0$
\end{enumerate}
If one of these conditions is fulfilled then $X\setminus P $ can be realized as an \'etale cover of the affine line of degree less than
or equal to $(2g-1)p$.
\end{theo}

\begin{proof} The equivalence of the conditions 1 and 2 is a direct consequence of Proposition~\ref{prop1}. A differential form $\omega\in\Omega(0)$ is nowhere vanishing on $X\setminus P $ if and only if it has a unique zero at $P$ of order $2g-2$. This easily follows from the identity $\deg(\omega)=2g-2$. In this case, the Riemann-Roch theorem asserts that the $k$-vector space $\Omega((2g-2)P)$ is one-dimensional, so that $\Ca(\Omega((2g-2)P))=0$ if and only if $\omega$ is exact, which shows the equivalence between the conditions 2 and 3. Let finally $x\in K$ be a rational function inducing an \'etale cover $X\setminus P \to\A$. It has a unique pole at $P$ and the last part of the theorem follows from Lemma~\ref{lemm1}, since in this case  we have 
$$\deg(x)=-v_P(x)\leq-\min\{(1-2g)p,v_P(\omega)+1\}=(2g-1)p$$
\end{proof}

\begin{rema}\label{rem2} The condition $\Omega((2g-2)P)\neq0$, which is equivalent to the identity $i((2g-2)P)=1$, is purely geometrical and can be restated as follows: let $J/k$ be the Jacobian variety of $X$. Fix a canonical divisor ${\mathcal K}_X$ on $X$. There is a natural morphism 
$$X\stackrel\varphi\longrightarrow J$$
 which sends a point $P\in X $ to the linear equivalence class defined by the divisor $(2g-2)P-{\mathcal K}_X$. Since any two canonical divisors on $X$ are linearly equivalent, the morphism $\varphi$ does not depend on the choice of ${\mathcal K}_X$. We then have the identity $i((2g-2)P)=1$ if and only if $\varphi(P)=0_J$, the origin of $J$. Remark that the curve $\varphi(X)$ may be singular, the singular points corresponding to the $(2g-2)$-torsion packets on $X$. Clearly, not any curve possesses a differential form having a unique zero. More precisely, let $\MM_g$ (resp. $\MM_{g,1}$) denote the coarse moduli space parametrizing isomorphism classes of smooth curves (resp. smooth $1$-pointed curves) of genus $g$ over $k$, viewed as a $k$-variety. Suppose that $g\geq2$. The {\it Weierstra{\ss} locus} $\W_g\subset\MM_{g,1}$ is the subset consisting of isomorphism classes of pairs $(X,P)$ where $P$ is a Weierstra{\ss} point of $X$, i.e. such that $i(gP)>0$. It is a divisor of $\MM_{g,1}$, so its dimension is $\dim\MM_{g,1}-1=3g-3$. If we restrict the canonical forgetfull morphism $\MM_{g,1}\to\MM_g$ to $\W_g$ we obtain a finite surjective morphism $\W_g\to\MM_g$. Let now $\Z_g\subset\W_g$ be the locus of pairs $(X,P)$ for which $i((2g-2)P)=1$. This is a closed subset of dimension $2g-1$. Its image in $\MM_g$ is also closed (the morphism $\MM_{g,1}\to\MM_g$ being proper) and has the same dimension. In general, $\Z_g$ has many irreducible components. One of them is the {\it hyperelliptic locus} $\H_g$, i.e. the closed subset consisting of isomorphisms classes of pairs $(X,P)$, where $X$ is an hyperelliptic curve of genus $g$ and $P$ is any Weierstra{\ss} point on it. See Remark~\ref{rem3} and \S5.2 for further investigations on $\Z_g$ and $\H_g$.
\end{rema}

This last result has many direct consequences. Recall that the curve $X$ is {\it ordinary} if the map $\Ca:\Omega(0)\to\Omega(0)$ is surjective (and thus, also injective). This can be restated by saying that the $p$-torsion subgroup $J[p]$ of the Jacobian variety of $X$ has order $p^g$. In general, $J[p]$ has order $p^{\sigma(X)}$, where the integer $\sigma(X)\in\{0,\dots,g\}$ is called the {\it $p$-rank} of the curve and there is a canonical bijection between $J[p]$ and the set of regular differential forms which are fixed by the Cartier operator (cf.~\cite{Serre_2},~\cite{Yui}).

\begin{coro}\label{non-ordinary} If $X\setminus P $ can be realized as an \'etale cover of the affine line then $p$ does not divide $2g-1$ and $X$ is not ordinary.
\end{coro}

\begin{proof} If $\omega$ is an exact differential form on $X$ then, for any $Q\in X $, the valuation of $\omega$ at $Q$ is an integer which is never congruent to $-1$ modulo $p$. In the present case, there exists an exact differential form $\omega\in\Omega(0)$ such that ${\rm div}(\omega)=(2g-2)P$ and thus $p$ does not divide $2g-1$. Moreover, the identity $\Ca(\omega)=0$ implies that the map $\Ca:\Omega(0)\to\Omega(0)$ is not injective, i.e. that the curve $X$ is not ordinary.
\end{proof}

\begin{coro}\label{(2g-2)-torsion} Suppose that $X$ has genus $g\geq2$. Then there exist finitely many elements $P\in X $ such that $X\setminus P $ can be realized as an \'etale cover of the affine line.\end{coro}

\begin{proof} Let $S$ be the subset of $X $ consisting of elements $P$ such that there exists an \'etale cover $X\setminus P \to\A$. If $S=\emptyset$ then there is nothing to prove. Suppose that $S$ is non-empty and fix an element $P_0$ in it. As above, let $J$ be the Jacobian variety of $X$. The map $S\to J$ which sends $P\in S$ to the linear equivalence class of the divisor $P-P_0$ is injective (it is just the restriction to $S$ of the albanese embedding $X\to J$ based at $P_0$). Moreover, the image of $S$ belongs to the $(2g-2)$-torsion subgroup of $J$. Indeed, if $P_1\neq P_0$ belongs to $S$ then, by the second condition of Theorem~\ref{th1}, there exist two regular differential forms $\omega_0$ and $\omega_1$ such that $\Div(\omega_i)=(2g-2)P_i$ for $i=0,1$. In particular, we have the relation $\omega_1=x\omega_0$ with $x\in K$ satisfying $\Div(x)=(2g-2)(P_1-P_0)$ and thus the divisor $(2g-2)(P_1-P_0)$ is principal, i.e. the order in $J$ of the image of the divisor $P_1-P_0$ divides $2g-2$, as claimed. The corollary is then immediate, since the $(2g-2)$-torsion subgroup of $J$ is finite.
\end{proof}

As we have seen, we can quite easily describe the curves  wich can be realized as an \'etale cover $X\setminus P \rightarrow\A$. The main inconvenient is that we don't have any control of its degree. In the rest of the paper, we develop a machinery which allows to overcome this problem. For the moment, we can give a first characterization: recall that the {\it gap sequence} $G(P)$ of $P$ is defined as
$$G(P)=\{n\in\bold N\,\,|\,\,L(nP)=L((n-1)P)\}$$
In other words, an integer $n$ belongs to $G(P)$ if and only if there does not exist any rational function on $X$ having a unique pole at $P$ of order $n$. This can be restated by saying that there exists a regular differential form having a zero at $P$ of order $n-1$. It is well known that $G(P)=\{n_1,\dots,n_g\}$ with $n_1=1$ (as soon as $g>0$) and $n_g<2g$. The condition $i((2g-2)P)=1$ can then be restated as $2g-1\in G(P)$. In the following we assume that there exists an \'etale cover $X\setminus P\to\A$. We say that a positive integer $m$ is an {\it admissible degree} if there exists such a cover of degree $m$.

\begin{prop}\label{minimal} With the above assumptions and notation, there exists an integer $n_0\in G(P)$ such that, for any positive integer $m$, the following conditions are equivalent:
\begin{enumerate}
\item $m$ is an admissible degree.
\item $m=np$ with $n=n_0$ or $n>n_0$ and $n\notin G(P)$.
\end{enumerate}
\end{prop}

\begin{proof} The existence of an \'etale cover $X\setminus P\to\A$ of degree $m$ is equivalent to the existence of a rational function $x$ having a unique pole at $P$ of order $m$ such that $dx$ has a unique zero at $P$ of order $2g-2$. We clearly have $m=np$ since otherwise $dx$ would have a pole at $P$ of order $m+1>0$, which is absurd. Suppose that $n$ does not belong to the gap sequence of $P$. We use the trick of Lemma~\ref{lemm1}: there exists a rational function $y$ having a unique pole at $P$ of order $n$ and we can choose a constant $c\in k^*$ such that $x_1=x-cy^p$ has a pole at $P$ of order strictly less than $m$. By construction, $dx_1=dx$ has a unique zero at $P$ and thus the rational function $x_1$ induces an \'etale cover $X\setminus P\to\A$. If we iterate this procedure we finally obtain a rational function $x_k\in L(n_0pP)$ with $n_0\in G(P)$. This shows that the minimal degree is equal to $p$ times an element of $G(P)$. Moreover, if $n>n_0$ is any integer not belonging to the gap sequence of $P$ then there exists a rational function $z$ having a unique pole at $P$ of order $n$ and the function $x_n+z^p$ induces an \'etale cover of the affine line of degree $d=np$. Suppose now that there exists two such covers induced by two rational functions $x_1,x_2$ of respective degrees $d_1=n_1p$ and $d_2=n_2p$ with $n_1,n_2\in G(P)$ and $n_1<n_2$. Since $i((2g-2)P)=1$ there exists a constant $c\in k^*$ such that $dx_1-cdx_2=0$. This can be restated by saying that the difference $x_1-cx_2=t^p$ is a $p$-power in $K$. But this is impossible, since the rational function $t$ would then have a unique pole at $P$ of order $n_2\in G(P)$.
\end{proof}

\begin{rema}\label{rem3} As it was pointed out in Remark~\ref{rem2}, the closed subset $\Z_g$ splits in many connected components. In fact, these components are completely determined by their gap sequence: two pointed curves $(X,P)$ and $(X',P')$ define two elements of $\Z_g$ belonging to the same connected component if and only if $G(P)=G(P')$. In particular, the locus $\E_g\subset\MM_{g,1}$ consisting of isomorphism classes of pairs $(X,P)$ such that $X\setminus P$ can be realized as an \'etale cover of the affine line will generally split in many connected components. For this reason, when studying $\E_g$ one may first of all fix the gap sequence. 
\end{rema}

\section{A sheaf interlude. Duality}

In this section we give some complements concerning the action of the Cartier operator on the set of meromorphic differential forms on $X$. For any integer $m$, set 
$$\iota(m)=\max\{n\in\ZZ\,\,|\,\,np\leq m \}$$
In other words, $\iota(m)$ is the integral part of the rational number $m/p$. We can then define a map
$$\aligned
{\DDiv}(X)&\stackrel\iota\longrightarrow{\DDiv}(X)\\
\sum_{P\in X}n_PP&\longmapsto\sum_{P\in X}\iota(n_P)P
\endaligned$$
For any divisor $D\in\DDiv(X)$, we have a homomorphism of $\OO$-modules 
$$F_*\Omega^1_X(D)\stackrel\Ca\longrightarrow\Omega^1_X(\iota(D))$$
which induces a $F^{-1}$-linear map $\mathscr C:\Omega(D)\to\Omega(\iota(D))$. Denote by $\ker(\Ca_D)$ its kernel and by $\coker(\Ca_D)$ its cokernel. Our first step is to give a cohomological interpretation of these $k$-vector spaces. Consider the sheaf of $\OO$-modules $\M(D)$ defined by the short exact sequence
\begin{equation}\label{short_1}
0\rightarrow\oX(\iota(D))\to F_*\oX(D)\rightarrow\M(D)\rightarrow 0
\end{equation}
For $D=0$, we recover the sheaf of {\it locally exact differential forms} on $X$ (cf.~\cite{Raynaud} and \cite{Madore}, where everything is formulated in terms of the relative Frobenius). In general $\M(D)$ is a locally free sheaf of $\OO$-modules of rank $p-1$. Remark moreover that $\M(pD)=\M(0)\otimes \oX(D)$. In order to state the next result, we need some more notation: consider the integer $\delta(D)$ defined by the relation
$$\delta(D)=l(D)-l(\iota(D))$$
The map $F:L(\iota(D))\to L(D)$ being injective, we obtain $\delta(D)\geq0$. We construct a second map from the set of divisors on $X$ into itself: let $\sigma:\bold Z\to\bold Z$ be the bijection defined by
$$\left\{\aligned
&\sigma(n)=-n\quad{\rm if}\; p|n\\
&\sigma(n)=-n-2\quad{\rm if}\; p|n+1\\
&\sigma(n)=-n-1\quad{\rm otherwise}
\endaligned\right.$$
Remark that we have the relations $\sigma^2=\id_\bold Z$ and $\sigma\circ\iota=\iota\circ\sigma$. As before, we obtain a map
$$\aligned
{\DDiv}(X)&\stackrel\sigma\longrightarrow{\DDiv}(X)\\
\sum_{P\in X}n_PP&\longmapsto\sum_{P\in X}\sigma(n_P)P
\endaligned$$

\begin{prop}\label{duality} The notation being as above, for any divisor $D$ on $X$, the following properties hold:
\begin{enumerate}
\item The $k$-vector spaces  $\ker(\Ca_D)$, $H^0(\M(\sigma(D)))$
and $H^1(\M(D))^\vee$ are canonically isomorphic.
\item There is a canonical injective map $\coker(\Ca_D)\hookrightarrow H^1(\M(\sigma(D)))$ which is an isomorphism if and only if $\delta(D)=0$. 
\end{enumerate}
\end{prop}

\begin{proof} There is a second fundamental short exact sequence of $\OO$-modules on which the relation between the Cartier operator and the sheaves $\M$ clearly appears, namely 
\begin{equation}\label{short_2}
0\rightarrow\M(\sigma(D))\stackrel{F_*d}\longrightarrow F_*\Omega^1_X(D)\stackrel\Ca\longrightarrow\Omega^1_X(\iota(D))\rightarrow0
\end{equation}
The associated long exact cohomology sequence splits in two parts,
$$0\to H^0(\M(\sigma(D)))\stackrel{F_*d}\longrightarrow\Omega(D)\stackrel\Ca\longrightarrow\Omega(\iota(D))\to\coker(\Ca_D)\to0$$
and 
$$0\to\coker(\Ca_D)\to H^1(\M(\sigma(D)))\stackrel{F^*d}\longrightarrow H^1(\Omega^1_X(D))\stackrel\Ca\longrightarrow H^1(\Omega^1_X(\iota(D)))\to0$$
from which we deduce the isomorphism $\ker(\Ca_D)\cong H^0(\M(\sigma(D)))$ and the injective map $\coker(\Ca_D)\hookrightarrow H^1(\M(\sigma(D)))$. Everything is clearly canonical. The above morphism 
$$H^1(\Omega^1_X(D))\stackrel\Ca\longrightarrow H^1(\Omega^1_X(\iota(D)))$$
is injective if and only if $\coker(\Ca_D)$ is isomorphic to $H^1(\M(\sigma(D)))$. Being the transpose of the Frobenius map
$$L(\iota(D))\stackrel F\longrightarrow L(D)$$
its injectivity is equivalent to the condition $L(\iota(D))^p=L(D)$, which can be restated as $\delta(D)=0$. Finally, the short exact sequence~\ref{short_1} gives the (truncated) long exact sequence in cohomology
$$H^1(\oX(\iota(D)))\to H^1(\oX(D))\to H^1(\M(D))\to 0$$
Its transpose is then given by
$$0\to H^1(\M(D))^\vee\to\Omega(D)\stackrel\Ca\longrightarrow\Omega(\iota(D))$$
so that $\ker(\Ca_D)$ is canonically isomorphic to $H^1(\M(D))^\vee$, which concludes the proof of the proposition.
\end{proof}

\begin{rema} There is a Serre duality between the sheaves $\M(D)$ and $\M(\sigma(D))$ which can be explicitely given as follows (cf.~\cite{Raynaud} for $D=0$): one easily checks that the map
$$\aligned \M(D)\times\M(\sigma&(D))\to\Omega^1_X\\
(x,y)&\mapsto\Ca(xdy)
\endaligned$$
is a perfect pairing of $\OO$-modules and thus it induces an isomorphism between $\M(\sigma(D))$ and $\Omega^1_X\otimes \M(D)^\vee$.
\end{rema}

\begin{coro}\label{duality_2} If $\delta(D)=0$ then $\coker(\Ca_D)$ and $\ker(\Ca_{\sigma(D)})^\vee$ are canonically isomorphic.
\end{coro}

\begin{proof} By the property 2 of Proposition~\ref{duality}, $\coker(\Ca_D)$ is canonically isomorphic to $H^1(\M(\sigma(D)))$. But the property 1 asserts that this last $k$-vector space is canonically isomorphic to $\ker(\Ca_{\sigma(D)})^\vee$.
\end{proof}

\begin{coro}\label{inequality} For any divisor $D$ on $X$, we have the inequality
$$\dim_k{\rm coker}(\Ca_{D})\leq i(\sigma(D))$$
\end{coro}

\begin{proof} We know from the property 2 of Proposition~\ref{duality} that there is an injective map $\coker(\Ca_D)\hookrightarrow H^1(\M(\sigma(D)))$. Moreover, by the property 1, the dual of this last vector space is canonically isomorphic to $\ker(\Ca_{\sigma(D)})$. We then obtain
$$\dim_k\coker(\Ca_D)\leq\dim_k H^1(\M(\sigma(D)))=\dim_k\ker(\Ca_{\sigma(D)})\leq i(\sigma(D))$$
\end{proof}

We now restrict to the case where $D=mP$ is based at a unique point. In general, the Riemann-Roch theorem asserts that if $\deg(D)$ is greater than or equal to $2g-1$ then $i(D)=0$. In particular, if $p$ does not divide $2g-1$ and if $m\leq-2g$ then $\sigma(mP)\geq2g-1$. Under these hypothesis, the above corollary implies that the $k$-vector space $\coker(\Ca_{mP})$ is trivial, i.e. the map 
$$\Omega(mP)\stackrel\Ca\longrightarrow\Omega(\iota(m)P)$$
is surjective. The following result describes the case $m=1-2g$ and generalizes Theorem~\ref{th1}.

\begin{prop}\label{th2_bis} For any point $P\in X $, the following conditions are equivalent:
\begin{enumerate}
\item The curve $X\setminus P $ can be realized as an \'etale cover of the affine line.
\item The $k$-vector space $\ker(\Ca_{(2g-2)P})$ is non-trivial.
\item The $k$-vector space $\coker(\Ca_{(1-2g)P})$ is non-trivial.
\end{enumerate}
\end{prop}

\begin{proof} The equivalence of the first two conditions is immediate, the second being just a restatement of the condition 3 in Theorem~\ref{th1}. Since $\delta((1-2g)P)=0$, Corollary~\ref{duality_2} asserts that $\coker(\Ca_{(1-2g)P})$ is canonically isomorphic to the dual of $\ker(\Ca_{(2g-2)P})$. In particular they are both trivial or non-trivial, which shows the equivalence between the last two conditions.
\end{proof}

\section{Etale covers of the affine line. Minimal degree}

From now on, we assume that there exists an \'etale cover $X\setminus P \to\A$ and we study the question of its minimal degree. The results of the previous section allow us to completely solve this problem.

\begin{theo}\label{minimal_2} Suppose that the curve $X\setminus P$ can be realized as an \'etale cover of the affine line. Then its minimal degree is the least integer $m=np>0$ such that $\Omega(nP)$ is contained in $\im(\Ca_{(1-2g)P})$.
\end{theo}

\begin{proof} The existence of an \'etale cover $X\setminus P \to\A$ of degree $n$ is equivalent to the existence of a rational function $x\in K$ having a unique pole at $P$ of order $n$ such that $dx$ has a unique zero at $P$ of ordrer $2g-2$. Clearly, in this case $n=mp$ is a multiple of $p$, since the converse would imply that $dx$ has a pole at $P$ of order $n+1$. Suppose that $x$ defines a cover of minimal degree $m=np$. We know from Proposition~\ref{minimal} that $n=n_i$ belongs to the gap sequence $G(P)=\{n_1,\dots,n_g\}$ of $P$. More precisely, by using the trick of Lemma~\ref{lemm1}, we can assume that there exists a uniformizer $z$ at $P$ such that in $\mathscr O_{X,P}$ we have the identity
$$x=a_iz^{-n_ip}+a_{i-1}z^{-n_{i-1}p}+\dots+a_1z^{-p}+b_1z^p+\dots+b_uz^{up}+z^{2g-1}h$$
with $h\in\mathscr O_{X,P}^*$, $u=\iota(2g-1)$ and $a_i\neq0$. In this case $x$ cannot be a $p$-power since there are no rational functions having a unique pole at $P$ of order $n_i$. Let $R$ denote the group of {\it repartitions} (or {\it ad\`eles}) on $X$. We know from~\cite{Serre_1} that the cohomology group $H^1(\oX(D))$ is canonically isomorphic to $R/(R(D)+K)$. We take $D=-jP$ with $j>0$. If $r=\{r_Q\}$ is the repartition defined by $r_Q=1$ for $Q\neq P$ and $r_P=z$ then the elements $r^{-n_g},\dots,r^{-n_1},r,r^2,\dots,r^{j-1}$ form a basis of $R/(R(D)+K)$ and the existence of a rational function $x$ admitting the above expansion at $P$ is equivalent to the relation $t^p=0$ in $R/(R((1-2g)P)+K)$, where we have set
$$t=a_ir^{-n_i}+\dots+a_1r^{-1}+b_1r+\dots+b_ur^u$$
The duality between $R/(R((1-2g)P)+K)$ and $\Omega((1-2g)P)$ is explicitely obtained via the perfect pairing
$$\langle t,\omega\rangle=\res_P(t\omega)$$
In particular, the adjoint property of the Cartier operator asserts that $t^p=0$ if and only if $\res_P(t\omega)=0$ for any $\omega\in\im(\Ca_{(1-2g)P})$. We know from Theorem~\ref{th2_bis} that the $k$-vector space $\coker(\Ca_{(1-2g)P})$ is one-dimensional. In particular, $\im(\Ca_{(1-2g)P})$ is the sub-$k$-vector space of $\Omega(\iota(1-2g)P)$ defined by the unique condition $\res_P(t\omega)=0$. If $\omega\in\Omega(\iota(1-2g)P)$ has a zero at $P$ of order greater than $n_i-1$ then $\res_P(t\omega)=0$ and thus $\omega\in\im(\Ca_{(1-2g)P})$. If $\omega$ has a zero of order $n_i-1$ (it exists, since $n_i$ belongs to the gap sequence of $P$) then $\res_P(t\omega)=a_i\neq0$ and thus $\omega$ does not belong to $\im(\Ca_{(1-2g)P})$. This concludes the proof of the theorem.
\end{proof}

\begin{coro}\label{generic} The assumptions being as in Theorem~\ref{minimal_2}, suppose that the $p$-rank of $X$ is equal to $g-1$. Then the minimal degree of the cover $X\setminus P \to\A$ is equal to $(2g-1)p$.
\end{coro}

\begin{proof} We have the inequality $\sigma(X)\leq g-\dim_k \ker(\Ca_0)$, so that $\ker(\Ca_0)$ is $1$-dimensional and spanned by $\Omega((2g-2)P)$. The identity
$$\Omega(\iota(1-2g)P)=\im(\Ca_{(1-2g)P})+\Omega(0)$$
is easily checked and follows from Riemann-Roch Teorem. The existence of a differential form $\omega\in\Omega((1-2g)P)$ for which $\Ca(\omega)$ has a unique zero at $P$ would imply that $\coker(\Ca_{(1-2g)P})$ is trivial, which is impossible by Proposition~\ref{th2_bis}. In particular, $\Omega((2g-2)P)$ is not conatined in $\im(\Ca_{(1-2g)P})$ and we just have to apply Theorem~\ref{minimal_2}.
\end{proof}

\section{Some explicit examples}

\subsection{Elliptic curves}

We start this serie of examples with the result which was the original motivation of the paper. Recall that a non-ordinary elliptic curve over $k$ is called {\it supersingular}.
 
\begin{prop}\label{elliptic} Let $E/k$ be an elliptic curve and denote by $0_E\in E$ its origin. Then, the following conditions are equivalent:
\begin{enumerate}
\item The curve $E$ is supersingular. 
\item There exists an \'etale cover $E\setminus 0_E\rightarrow\A$.
\item There exists an \'etale cover $E\setminus 0_E\rightarrow\A$ of degree $p$.
\end{enumerate}
\end{prop}

\begin{proof} Since $(2g-1)p=p$, Theorem~\ref{th1} asserts that if there exists an \'etale cover $E\setminus 0_E\rightarrow\A$, then we can suppose that its degree $n$ is less than or equal to $p$. The condition $n<p$ would lead to a cover $E\rightarrow\P$ tamely ramified above only one point, which is impossible, so that we obtain $n=p$. This shows the equivalence of the conditions 2 and 3. In this special case, the space $\Omega(0)$ is one-dimensional and the third condition of Theorem~\ref{th1} can be rewritten as $\Ca(\Omega(0))=0$, i.e. the curve $E$ is supersingular. This shows the equivalence between the conditions 1 and 2 above.
\end{proof}

The following table performs the complete list of the supersingular elliptic curves and of the associated \'etale covers of the affine line for any prime $p\leq17$.

\vskip.5cm
\begin{center}
\begin{tabular}{|c|c|c|c|}
\hline
$p$ & Equation & $j$-invariant & Cover\\
\hline
$2$ & $y^2+y=x^3$ & $0$ & $x$\\
\hline
$3$ & $y^2=x^3-x$ & $0$ & $y$\\
\hline
$5$ & $y^2=x^3-1$ & $0$ & $xy$\\
\hline
$7$ & $y^2=x^3-x$ & $6$ & $(x^2+4)y$\\
\hline
$11$ & $y^2=x^3-1$ & $0$ & $(x^4+6x)y$\\
\hline
$11$ & $y^2=x^3-x$ & $1$ & $(x^4+6x^2+10)y$\\
\hline
$13$ & $y^2=x^3+x+4$ & $5$ & $(x^5+6x^3+11x^2+2x+3)y$\\
\hline
$17$ & $y^2=x^3-1$ & $0$ & $(x^7+9x^4+11x)y$\\
\hline
$17$ & $y^2=x^3+x-1$ & $8$ & $(x^7+8x^5+9x^4+11x^3+12x^2-x+2)y$\\
\hline
\end{tabular}
\end{center}
\vskip.5cm

\begin{coro}\label{cor4} If $E\setminus 0_E$ can be realized as an \'etale cover of the affine line then $E$ can be defined over $\bold F_{p^2}$.
\end{coro}

\begin{proof} We know from the general theory of elliptic curves that if $E$ is supersingular then its $j$-invariant $j(E)$ belongs to $\bold F_{p^2}$ (cf.~\cite{Silverman}). The corollary follows from Proposition~\ref{elliptic} and from the fact that $\bold F_p(j(E))$ is the minimal field of definition of $E$ (cf.~{\it loc. cit.}).
\end{proof}

\begin{coro}\label{cor5} Any elliptic curve $E/k$ can be realized as a degree $p$ cover of the projective line unramified outside at most two points and wildly ramified above them.\end{coro}

\begin{proof} Following Proposition~\ref{elliptic}, we can suppose that the curve $E$ is ordinary. In this case, it has $p$-torsion and we can choose $P\in E\setminus 0_E$ such that the divisor $D=p(P-0_E)$ is principal, say $\Div(x)=D$ with $x\in K$. Clearly, $x$ is not a $p$-power (because $P-0_E$ cannot be principal) and so the induced cover $E\to\P$ has degree $p$, is generically \'etale and wildly ramified above $0$ and $\infty$. The other ramification points correspond to the zeroes of $dx$. Remark that $D'=\Div(dx)\geq D$. Since $\deg(D')=\deg(D)=0$, we finally obtain $D=D'$ and thus there is no other ramification.
\end{proof}

In the supersingular case, the degree $p$ cover $E\to\P$ is canonical and can be defined over the field of moduli of the curve. This property is no longer true if $E$ is ordinary, since the cover depends on the choice of a $p$-torsion point. Using Proposition~\ref{duality}, we now show that there always exists a canonical degree $p$ cover $E\to\P$ unramified outside at most three points and wildly ramified above only one of them.

\begin{prop}\label{elliptic_general} Suppose that $p>2$. Then, any elliptic curve $E/k$ can be realized as a degree $p$ cover of the projective line unramified outside at most three points, wildly ramified above one of them and tamely ramified above the others.
\end{prop}

\begin{proof} We know from Corollary~\ref{inequality} that the map
$$\Omega(-20_E)\stackrel\Ca\longrightarrow\Omega(0)$$
is surjective. Since $\Omega(0)$ is $1$-dimensional and $\Omega(-20_E)$ is $2$-dimensional, we deduce that there exists a differential form $\omega\in\Omega(-20_E)$ such that $\Ca(\omega)=0$. Moreover, Lemma~\ref{lemm1} implies that there exists a rational function $x$ having a unique pole at $0_E$ of degree $n\leq p$ such that $dx=\omega$. The case $n<p$ is impossible, since in this case $dx$ would have a pole at $0_E$ of order $n+1>2$. We then obtain $n=p$. The zeroes of $\omega$ are the ramification points of the cover $E\to\P$ induced by $x$. Since the valuation of $\omega$ at $0_E$ (which is a ramified point) is greater than or equal to $-2$, we see that there are at most two other ramified points and that the valuation of $\omega$ at these points is less than or equal to $2$ (less than or equal to $1$ for $p=3$, since the valuation of $\omega$ at any point is never congruent to $-1$ mod $p$). In particular, the corresponding ramification indices are less than or equal to $3$ (less than or equal to $2$ for $p=3$) and thus the ramification is tame.
\end{proof}

The cover $E\to\P$ constructed in the proof of Proposition~\ref{elliptic_general} is canonical and can be defined over the field of moduli of the curve. We now show how to explicitly find its ramified points. Since we are assuming that $p>2$, the curve $E$ can be given by an affine Weierstra{\ss} equation
$$y^2=x(x-1)(x-\lambda)$$
the identity element $0_E$ corresponding to the point at infinity. The $k$-vector space $\Omega(-20_E)$ is spanned by the differential forms $\omega_1=dx/y$ and $\omega_2=xdx/y$. Setting $m=\frac{p-1}2$ and
$$c_n=(-1)^n\sum_{i=0}^n\binom mi\binom m{n-i}\lambda^{m-i}$$
for any $n\in\{0,\dots,m\}$, we then have the identities
$$\left\{\aligned
&\Ca(\omega_1)=F^{-1}(c_m)\omega_1\\
&\Ca(\omega_2)=F^{-1}(c_{m-1})\omega_1
\endaligned\right.$$
where $F$ is the Frobenius map. Since $\Ca_{-20_E}$ is surjective (cf. Proposition~\ref{inequality}), we see that $c_m$ and $c_{m-1}$ have no common roots (as polynomials on $\lambda$). A direct computation shows that this fact is far from being trivial. The $k$-vector space $\ker(\Omega(-20_E))$ is spanned by the differential form
$$\omega=c_{m-1}\omega_1-c_m\omega_0=(c_mx-c_{m-1})\omega_1$$
and its zeroes are the ramified points of the cover $E\to\P$ which are different from $0_E$. In particular, we see that if $c_m=0$ then the curve is supersingular (this is a well-know criterion, cf.~\cite{Silverman}) and there are no ramified points different from $0_E$. This is just a restatement of Proposition~\ref{elliptic}. Suppose that $c_m\neq0$ and set $c=c_m^{-1}c_{m-1}$. The differential form $\omega$ has a unique zero $P$ if and only if $c\in\{0,1,\lambda\}$. The corresponding degree $p$ cover $E\to\P$ is unramified outside two points. In the remaining cases, $\omega$ has two zeroes and the cover is unramified outside two or three points.

\subsection{hyperelliptic curves in characteristic $p>2$}

As we have seen in Remarks~\ref{rem2} and~\ref{rem3}, from which we keep the notation, the locus $\Z_g\subset\MM_{g,1}$ splits in many irreducinble components. The aim of this section is to study more closely its connected component $\H_g$, the hyperelliptic locus. Remark that if $\H_g'$ denotes the image of $\H_g$ in $\MM_g$ by the (restriction of the) forgetfull map then the morphism $\H_g\to\H_g'$ is generically \'etale of degree $2g+2$. We assume that the characteristic $p$ of $k$ is strictly greater than $2$. Let $X$ be an hyperelliptic curve. In this case, we have $i((2g-2)P)>0$ if and only if $P$ is a Weierstra{\ss} point. From now on, we fix such a point $P$ (there are $2g+2$ of them). We know from the general theory (cf.~\cite{Yui}) that there exists a monic polynomial $f(x)\in k[x]$ of degree $2g+1$ having pairwise distinct roots such that the equation
$$y^2=f(x)$$
defines a non-singular (affine) model of $X\setminus P $ in $\bold A^2_k$. Remark that for $g\geq2$ the closure of this model in $\bold P^2_k$ is singular. Clearly, $f(x)$ is not unique but any other such polynomial $f_1(x)$ can be expressed as 
$$f_1(x)=a^{-2g-1}f(ax+b)$$ with $a\in k^*$ and $b\in k$. If $p$ does not divide $2g+1$ then we can moreover suppose that the trace of $f(x)$ (i.e. the coefficient of $x^{2g}$) is zero. In this case, $f(x)$ is uniquely determined up to a transformation of the type $f(x)\mapsto\lambda^{-2g-1}f(\lambda x)$ with $\lambda\in k^*$. If we set $f(x)=a_0+\dots+a_{2g-1}x^{2g-1}+x^{2g+1}$ then we can associate to it the $2g$-ple $(a_0,\dots,a_{2g-1})\in\bold A^n(k)$. Two such polynomials $f(x)$ and $g(x)$, corresponding respectively to $(a_0,\dots,a_{2g-1})$ and $(b_0,\dots,b_{2g-1})$, determine two isomorphic $1$-pointed curves if and only if there exists an element $\lambda\in k^*$ such that
$$(b_0,\dots,b_{2g-1})=(\lambda^{2g+1}b_0,\lambda^{2g}b_1,\dots,\lambda^2b_{2g-1})$$
This construction shows that $\H_g$ can naturally be identified with an open subset of the weighted projective space $\bold P(2,3,\dots,2g+1)$ (which is the complement of the divisor corresponding to the polynomials with zero discriminant).

\begin{prop}\label{hyperelliptic} Let $X\setminus P $ be a $1$-pointed hyperelliptic curve of genus $g$ defined by an affine Weierstra{\ss} equation $y^2=f(x)$. Consider the polynomial $h(x)\in k[x]$ defined by
$$h(x)=f(x)^{p-1\over2}=\sum_{n\geq0}c_nx^n$$
The following conditions are equivalent:

\begin{enumerate}
\item The curve $X\setminus P $ can be realized as an \'etale cover of the affine line.
\item $c_{p-1}=\dots=c_{gp-1}=0$.
\item There exists a polynomial $H(x)\in k[x]$ with $H'(x)=h(x)$.
\end{enumerate}

\end{prop}

\begin{proof} For any $n\in\{1,\dots,g\}$, set
$$\omega_n=x^{n-1}\frac{dx}y$$
These differential forms are everywhere regular on $X$ and define a (canonical) basis of the $k$-vector space $\Omega(0)$.
We moreover have the identity $v_P(\omega_n)=2(g-n)$. In particular, $\omega_1$ has a unique zero at $P$ of order $2g-2$ and spans $\Omega((2g-2)P)$. 
Following Theorem~\ref{th1}, the curve $X\setminus P $ can be realized as an \'etale cover of the affine line if and only if $\Ca(\omega_1)=0$. An easy and explicit computation (see for example~\cite{Yui}) gives
$$\Ca(\omega_1)=y^{-1}\Ca(h(x)dx)=\sum_{n=1}^gF^{-1}(c_{np-1})\omega_n$$
where as usual $F$ denotes the Frobenius map. In particular, since $F$ defines an isomorphism of $k$, we have $\Ca(\omega_1)=0$ if and only if $c_{p-1}=\dots=c_{gp-1}=0$. Moreover, the existence of the polynomial $H(x)$ is equivalent to the identity $\Ca(h(x)dx)=0$ and the proposition follows.
\end{proof}

\begin{coro}\label{moduli} The locus $\E_g^\hyp=\E_g\cap\H_g$ consisting of isomorphism classes of pairs $(X,P)$ such that $X$ is hyperelliptic and $X\setminus P$ can be realized as an \'etale cover of the affine line is the intersection of $g$ divisors of $\H_g$. In particular, if it is non empty then its dimension is greater than or equal to $g-1$.
\end{coro}

\begin{proof} We assume that $p$ does not divide $2g+1$, this last case can be treated similarly. As we have just seen, we can identify $\H_g$ with an open subset of the weighted projective space $\bold P(2,3,\dots,2g+1)=\proj(k[a_0,\dots,a_{2g-1}])$, where $a_i$ is homogeneous of degree $2g+1-i$. In this case, Proposition~\ref{hyperelliptic} asserts that the locus $\E_g^\hyp$ is defined by the $g$ equations $c_{p-1}=\dots=c_{gp-1}=0$. Now, $c_{ip-1}\in k[a_0,\dots,a_{2g-1}]$ is homogeneous and the equation $c_{ip-1}=0$ defines a divisor in $\H_g$. The corollary is then immediate.
\end{proof}

\begin{coro}\label{p=3} Suppose that $k$ has characteristic $3$ and that $g=3n$. Then $\E_g^\hyp$ can be naturally be identified with an open subset of the weighted projective space $\bold P(3,4,6,7,\dots,3i,3i+1,\dots,6n,6n+1)$. In particular, it is connected of dimension $4n-1$.
\end{coro}

\begin{proof} The condition $g=3n$ is equivalent to $3\not|(2g-1)(2g+1)$. We know from Corollary~\ref{non-ordinary} that the condition $3\not|2g-1$ is necessary for $\E_g$ not being empty while $3\not|2g+1$ implies that $\H_g$ can be naturally identified with an open subset of the weighted projective space $\bold P(2,3,\dots,2g+1)$. Following the notation of Proposition~\ref{hyperelliptic}, we find $h(x)=f(x)=a_0+a_1x+\dots+a_{2g-1}x^{2g-1}+x^{2g+1}$, so that $\E_g^\hyp$ is defined by the equations $a_2=a_5=\dots=a_{2g-1}=0$. But this is just the weighted projective space $\bold P(3,4,6,7,\dots,2g,2g+1)\subset\bold P(2,3,\dots,2g+1)$.
\end{proof}

\begin{rema} This last result shows that in general the dimension of $\E_g^\hyp$ can be strictly greater than the bound $g-1$ given in Corollary~\ref{moduli}. Nevertheless, if the characteristic of $k$ is sufficiently large with respect to $g$ then we actually find a dimension of $g-1$.
\end{rema}

The last condition of Proposition~\ref{hyperelliptic} is particularly adapted for explicit calculations. From now on, we always assume that the the degree of the polynomial $H(x)$ is as small as possible, that is $\deg(H)=\deg(h)+1=(2g+1)\frac{p-1}2+1$.

\begin{prop}\label{explicit_cover} Suppose that the conditions of Proposition~\ref{hyperelliptic} are fulfilled and consider the polynomial  
$$H_1(x)=H(x)-\sum_{f(u)=0}H(u)\prod_{f(v)=0, v\neq u}\left({x-v\over u-v}\right)^p$$
Then the rational function $t=y^{-p}H_1(x)$ induces an \'etale cover $X\setminus P \to\A$ of minimal degree.
\end{prop}

\begin{proof} Setting $z=y^{-p}H(x)$, we find
$$dz=y^{-p}dH(x)=y^{-p}h(x)dx=y^{-p}f(x)^{p-1\over2}dx=y^{-p}y^{p-1}dx=y^{-1}dx=\omega_1$$
In general, the rational function $z$ has more than a unique pole at $P$, and thus it {\it does not} induce an \'etale cover $X\setminus P \to\A$. Remark that the rational function $t-z$ belongs to $K^p$, and thus $dt=dz=\omega_1$. Let $Q\in X(k)$ different from $P$. If $y(Q)\neq0$ then $t$ is regular at $Q$. Suppose that $y(Q)=0$, so that $Q$ is a Weierstra{\ss} point and corresponds to a root $u=x(Q)$ of the polynomial $f(x)$. In this case we obtain $v_Q(y)=1$ and $v_Q(x-u)=2$. By construction, we have $H_1(u)=0$ and $u$ is a root of $H_1'(x)=h(x)$ of multiplicity $\frac{p-1}2$. In particular, we obtain the factorization $H_1(x)=(x-u)^{\frac{p+1}2}H_2(x)$ with $H_2(u)\neq0$, from which we deduce the identity
$$v_Q(t)=v_Q(y^{-p}H_1(x))=\frac{p+1}2v_Q(x-u)-pv_Q(y)=1$$
This shows that $t$ is everywhere regular on $X\setminus P $ and thus, since $dt=\omega_1$, it induces an \'etale cover $X\setminus P \to\A$. One easily checks that the degree of the polynomial $H_1(x)$ is less than or equal to $2gp$. Since $v_P(x)=-2$ and $v_P(y)=-2g-1$, we finally obtain
$$\aligned 
n=&\deg(t)=-v_P(t)=pv_P(y)-v_P(H_1(x))=\\
=&2\deg(H_1)-(2g+1)p\leq(2g-1)p
\endaligned$$
Remark that the above identities imply that $n$ is an odd integer. Suppose that there exists an \'etale cover $X\setminus P \to\A$ of degree $m<n$. It corresponds to a rational function $v\in K$ having a unique pole at $P$ of order $m$ and satisfying $dv=c\omega_1$ with $c\in k^*$. We can clearly assume that $c=1$, so that $dt=dv=\omega_1$. But in this case we obtain $t-v=s^p$ where $s\in K$ has a unique pole at $P$ of odd order $p^{-1}n\leq2g-1$. This is impossible, since the gap sequence of $P$ is $G(P)=\{1,3,\dots,2g-1\}$. 
\end{proof}

\begin{coro}\label{generic-hyperelliptic} Suppose that there exists an \'etale cover $X\setminus P \to\A$. Then its minimal degree is strictly less than $(2g-1)p$ if and only if
$$\sum_{f(u)=0}{H(u)\over f'(u)^p}=0$$
\end{coro}

\begin{proof} We know from Proposition~\ref{explicit_cover} that the minimal degree of the cover is equal to $2\deg(H_1)-(2g+1)p$. In particular, it is strictly less than $(2g-1)p$ if and only if $\deg(H_1)<2gp$. The corollary follows from the fact that the above expression is just the coefficient of $x^{2gp}$ in $H_1(x)$.
\end{proof}

\subsection{Curves of genus two} 

We now illustrate these results in the case of genus two curves. They all are hyperelliptic and thus $\E_2=\E_2^\hyp$. Except from the case $p=5$, the moduli space $\H_2\subset\MM_{2,1}$ can be identified with an open subset of the weighted projective space $\bold P(2,3,4,5)$ and $\E_2$ is the intersection of $2$ divisors of $\H_2$. In the following, we say that $X$ is {\it supersingular} (resp. {\it supserspecial}) if the Cartier operator is nilpotent (resp. trivial) on $\Omega(0)$. Remark that $X$ is supersingular if and only if its Jacobian variety has no $p$-torsion, i.e. if $\sigma(X)=0$. Let $P\in X$ be one of the six Weierstra{\ss} points and suppose that there exists an \'etale cover $X\setminus P\to\A$. The gap sequence of $P$ being $G(P)=\{1,3\}$, Proposition~\ref{minimal} asserts that the minimal degree of the cover is equal either to $p$ or to $3p$. Moreover, we know from Corollary~\ref{non-ordinary} that $\sigma(X)\leq1$. If this last inequality is in fact an equality then Corollary~\ref{generic} implies that the minimal degree is equal to $3p$. If $X$ is not superspecial then the converse is also true. Indeed, we have the following result:

\begin{prop}\label{g=2} Let $X/k$ be a smooth non-superspecial curve of genus $2$ and suppose that $X\setminus P$ can be realized as an \'etale cover of the affine line. Then its minimal degree is equal to $p$ if and only if $X$ is supersingular.
\end{prop}

\begin{proof} We know from the above discussion that if $\sigma(X)=1$ then the minimal degree is $3p$. In particular, the supersingularity is necessary for the existence of a degree $p$ \'etale cover of the affine line. Suppose now that $\sigma(X)=0$ and set $\Omega(0)=\langle\omega_1,\omega_2\rangle$ with $\Div(\omega_1)=2P$. Remark that $v_P(\omega_2)=0$. Since the curve is non-superspecial and $\Ca(\omega_1)=0$ we deduce that $\Ca(\omega_2)\neq0$, say $\Ca(\omega_2)=a\omega_1+b\omega_2$ with $a\neq0$ or $b\neq0$. One easily checks that if $b\neq0$ then there exists a differential form $\omega\in\Omega(0)$ fixed by the Cartier operator and thus it corresponds to a non-trivial $p$-torsion point of the Jacobian of $X$. But this is impossible by the supersinguarity assumption. We then obtain $\Ca(\omega_2)=a\omega_1$ with $a\neq0$ so that $\omega_1\in\im(\Ca_0)$, which can be restated as $\Omega(2P)\subset\im(\Ca_0)$. Since $\im(\Ca_0)\subset\im(\Ca_{(1-2g)P})$, Theorem~\ref{minimal_2} asserts that the minimal degree is less than or equal to $2p$ and thus equal to $p$.
\end{proof}

The superspecial case is more delicate and for the moment we don't have any simple characterization. We can only affirm that in this case, for any Weierstra{\ss} point $P\in X$, there always exist a cover $X\setminus P\to\A$ but we cannot control its degree. A direct explicit computation shows that both $p$ and $3p$ can occur as minimal degrees (cf. Example~5.12). What seems reasonable is that for {\it at least} one of the Weierstra{\ss} points there exists such a cover of degree $p$. Remark finally that in genus two the curve $X$ is superspecial if and only if there exist two distinct points $P$ and $Q$ such that both $X\setminus P$ and $X\setminus Q$ can be realized as an \'etale cover of the affine line (since there exist two linearly independent regular differential forms on $X$ on which the Cartier operator acts trivially).

\begin{exem} We close this section and the whole paper with some explicit examples. First of all, there are no genus two curves $X$ in characteristic $3$ such that $X\setminus P$ can be realized as an \'etale cover of the affine line. Indeed, in this case $3$ divides $2g-1$ (cf. Corollary~\ref{non-ordinary}). In general, the locus $\E_2=\E_2^\hyp$ is one-dimensional, i.e. it is a curve.

For $p=5$ we find the following (uni-)versal family
$$X_a:y^2=x+ax^2+2a^2x^3+x^5$$
with $a\in k$ satisfying the inequality $a^4+3\neq0$, which ensures the non-singularity. We have $X_a\cong X_b$ if and only if there exists $\zeta\in k$ with $\zeta^4=1$ such that $b=\zeta a$. In particular, by setting $t=a^4$, we can identify $\E_2$ with ${\rm Spec}(k[t,1/(t+3)])$, which is isomorphic to the affine line minus one point. Since $1-2g=-3$ and $\iota(1-2g)=-1$ the Cartier operator maps $\Omega(-3P)$ to $\Omega(-P)=\Omega(0)$. The $k$-vector space $\Omega(0)$ is spanned by the differential forms $dx/y$ and $xdx/y$ and $\Omega(-3P)=\Omega(0)\oplus\langle x^2dx/y,dy\rangle$. The matrix associated to the map $\Ca_{-3P}$ with respect to this basis is given by
$$\left(
\begin{matrix}
0 & 2a & 1 & 0\\
0 & -a^2 & 2a & 0
\end{matrix}
\right)$$
For $a\neq0$ we easily see that $n=3$ is the least integer for which $\Omega(nP)$ is contained in $\im(\Ca_{-3P})$, so that Theorem~\ref{minimal_2} asserts that the minimal degree of the cover is $15$. Remark that in this case the $5$-rank of the curve $X_a$ is equal to $1$, as predicted by Proposition~\ref{g=2}. For $a=0$ we find the unique superspecial curve in characteristic $5$. The minimal degree is then equal to $5$ and the cover is Galois.

The case $p=7$ leads to the following (uni-)versal family
$$X_{a,b}=a^3+b^3x+abx^2+x^5$$
with $a,b\in k$ such that $a^4+3b^5\neq0$. We have $X_{a,b}\cong X_{c,d}$ if and only if $a^4d^5=c^4b^5$, i.e. if there exists $u\in k^*$ such that $c=u^5a$ and $d=u^4b$. One easily checks that $\E_2$ can then be identified with the weighted projective space $\bold P(4,5)$ minus one point, that is to the affine line. Concerning the minimal degree, by taking the same basis than in the previous example we find the matrix
$$\left(
\begin{matrix}
0 & 3a^2(b^5+a^4) & 3ab^2(b^5+a^4) & 0\\
0 & 3ab & 3b^3 & 0
\end{matrix}
\right)$$
If $ab\neq0$ then, as above, $n=3$ is the least integer for which $\Omega(nP)$ is contained in $\im(\Ca_{-3P})$ so that the minimal degree is equal to $21$. A direct computation shows that the curve $X_{a,b}$ is supersingular if and only if $ab=0$. For $b=0$ we can set $a=1$ and the corresponding curve is not superspecial. In this case, Proposition~\ref{g=2} asserts that the minimal degree is equal to $7$. For $a=0$ and $b=1$ we obtain a superspecial curve and the minimal degree is equal to $21$.

\end{exem}

\section*{\large{Acknowledgment}} I would like to thank the organizers and the participants of the GTEM midterm review in Leiden (June 2002) for giving me the opportunity to give a talk on this subject. At that moment, the question were only studied in the genus one situation and the later investigations have been positively motivated by many discussions during this meeting. I thank F. Pop for inviting me for in Bonn (July 2002 and July-August 2003); all the general results were obtained during these periods. I am particularly endebded with J. Stix for encouraging me since the very begining and for his interest and criticism in reading the preliminary versions of the paper.

\end{document}